\documentclass[12pt]{amsart}
\usepackage{amssymb,latexsym,eufrak,amsmath,amscd,graphicx}
 \usepackage[all]{xy}

\theoremstyle{plain}
\newtheorem{theorem}{Theorem}[section]

\newtheorem{corollary}[theorem]{Corollary}

\theoremstyle{definition}
\newtheorem{definition}[theorem]{Definition}

\newtheorem{remark}[theorem]{Remark}
\newtheorem{example}[theorem]{Example}

\newtheorem{conjecture}[theorem]{Conjecture}

\newcommand{\lra}{\longrightarrow}

\newcommand{\AAA}{\mathbb{A}}
\newcommand{\PP}{\mathbb{P}}

\newcommand{\RR}{\mathbb{R}}
\newcommand{\NN}{\mathbb{N}}
\newcommand{\ZZ}{\mathbb{Z}}
\newcommand{\CC}{\mathbb{C}}
\newcommand{\QQ}{\mathbb{Q}}

\newcommand{\OO}{\mathcal{O}}
\newcommand{\JJ}{\mathcal{J}}

\newcommand{\fra}{\frak{a}}

\newcommand{\frmm}{\frak{m}}
\newcommand{\frakm}{\frak{m}}

\newcommand{\Spec}{\textnormal{Spec}}

\newcommand{\ord}{\textnormal{ord}}
\newcommand{\codim}{\textnormal{codim}}
\newcommand{\mult}{\textnormal{mult}}
\newcommand{\pr}{\prime}

\newcommand{\Cont}{\textnormal{Cont}}

\newcommand{\val}{\textnormal{val}}

\newcommand{\lc}{\textnormal{lc}}
\newcommand{\Aut}{\textnormal{Aut}}
\newcommand{\Bir}{\textnormal{Bir}}

\begin{document}

\title{Invariants of Singularities of Pairs}

\author[L.Ein]{Lawrence Ein}
\address{Department of Mathematics \\ University
of Illinois at Chicago, \hfil\break\indent  851 South Morgan Street
(M/C 249)\\ Chicago, IL 60607-7045, USA \hfil\break\indent and \hfil
\break\indent  Department of Mathematics
\\ University of California at Irvine, \hfil\break\indent  Irvine,
CA 92697-3875, USA} \email{ein@math.uic.edu}

\author[M. Musta\c{t}\v{a}]{Mircea Musta\c{t}\v{a}}
\address{Department of Mathematics \\ University of Michigan \\
Ann Arbor, MI 48109,  USA}
\email{mmustata@umich.edu}

\thanks{Research of Ein and Musta\c{t}\v{a}
  was partially supported by the NSF under grants DMS  0200278
  and 0500127.}

\maketitle


\noindent {\bf Abstract.} Let $X$ be a smooth complex variety and
$Y$ be a closed subvariety of $X$, or more generally, a closed
subscheme of $X$. We are interested in invariants attached to the
singularities of the pair $(X, Y)$. We discuss various methods to
construct such invariants, coming from the theory of multiplier
ideals, D-modules, the geometry of the space of arcs and
characteristic $p$ techniques. We present several applications of
these invariants to algebra,  higher dimensional birational geometry
and to singularities.

\section{Introduction}
Let $X$ be a smooth complex variety of dimension $n$ and $Y$ be a
closed subvariety of $X$ (or more generally, a closed subscheme of
$X$). We are interested in studying the singularities of the pair
$(X, Y)$. The general setup is to assume only that $X$ is normal and
$\QQ$-Gorenstein, as in \cite{Kollar1}. However, several of the
approaches we will discuss become particularly transparent if we
assume, as we do, the smoothness of the ambient variety. The
following are some examples of pairs.

\begin{example}
\begin{enumerate}
\item[(i)] $X = \CC^n$ and $Y$ is a hypersurface defined by an equation
$\,\,\,\,\,\,\,\,$ $f(x_1, \dots, x_n) = 0$. For instance $f$ can be
the Fermat hypersurface $x_1^m + x_2^m + \dots + x_n^m = 0$, which
has an isolated singularity of multiplicity $m$ at the origin.

\item[(ii)] If
$X$ is a smooth projective variety and $L$ is a line bundle on $X$,
then we take  $Y$ to be the base locus of the complete linear system
$|L|$, i.e. $Y = \bigcap_{D \in |L|} D$.

\item[(iii)] Let $X$ be a smooth affine variety with coordinate ring $R$.
If $I \subseteq R$ is an ideal, then we take $Y$ to be the closed
subscheme defined by $I$.
\end{enumerate}
\end{example}

In what follows we present various invariants attached to such pairs
and we discuss some of their applications. Our main point is that
the same invariants that play an important role in higher
dimensional algebraic geometry arise also in several other
approaches to singularities.

\section{Multiplier ideals}

Multiplier ideals were first introduced by J. Kohn for solving
certain partial differential equations. Siu and Nadel introduced
them to complex geometry.  We discuss below these ideals in the
context of algebraic geometry.

Let $X$ be a smooth complex affine variety and $Y$ be a closed
subscheme of $X$. Suppose that the ideal of $Y$ is generated by
$f_1,\ldots,f_m$, and let $\lambda$ be a positive real number. We
define the multiplier ideal of $(X, Y)$ of coefficient $\lambda$ as
follows:

$$ \JJ(X, \lambda\cdot Y) = \left\{g \in \OO_X\mid
\frac{|g|^2}{(\sum_{i=1}^{m} |f_i|^2)^{\lambda}}\,{\rm is}\,{\rm
locally}\,{\rm integrable}\right\}.
$$

\begin{example}\label{SNC}
Let $X = \CC^n$ and let $Y$ be the closed subscheme of $X$ defined
by $f=x_1^{a_1}\cdots x_n^{a_n}$. Then
$$\JJ(X, \lambda \cdot Y) =
(x_1^{\lfloor\lambda a_1\rfloor}\dots x_n^{\lfloor\lambda
a_n\rfloor}),$$ where $\lfloor\alpha\rfloor$ denotes the integer
part of $\alpha$. Equivalently, if $H_i$ is the divisor defined by
$x_i=0$, then $g$ is in $\JJ(X,\lambda\cdot Y)$ if and only if
$$ \ord_{H_i} g \ge \lfloor\lambda a_i\rfloor, $$
for $i= 1, \dots, n$.
\end{example}

We can use a log resolution of singularities and the above example
to give in general a more geometric description of the multiplier
ideals of $(X, Y)$. By Hironaka's Theorem there is a log resolution
of singularities of the pair $(X,Y)$, i.e. a proper birational
morphism
$$ \mu : X' \lra X$$
with the following properties. The variety $X'$ is smooth,
$\mu^{-1}(Y)$ is a divisor, and the union of $\mu^{-1}(Y)$ and the
exceptional locus of $\mu$ has simple normal crossings. The relative
canonical divisor $K_{X'/X}$ is locally defined by the determinant
of the Jacobian $J(\mu)$ of $\mu$, hence it is supported on the
exceptional locus of $\mu$.
 We write
$\mu^{-1}(Y)=
 \sum_{i=1}^{N} a_i E_i$ and
 $K_{X'/X} = \sum_{i=1}^{N}k_i E_i$, where the $E_i$ are
distinct smooth irreducible divisors in $X'$ such that
$\sum_{i=1}^{N} E_i$ has only simple normal crossing singularities.

Suppose that $x_1,\dots, x_n$ are local coordinates in $X$ and $y_1,
y_2, \dots, y_n$ are local coordinates for an open set in $X'$. Note
that
\begin{equation}
\mu^*dx_1\cdots dx_nd\overline{x}_1\cdots d\overline{x}_n =
|det(J(\mu)|^2 dy_1\cdots dy_nd\overline{y}_1 \cdots
d\overline{y}_n.
\end{equation}
The local integrability of a function $g$ on $X$ can be expressed as
a local integrability condition on $X'$ via the change of variable
formula. This reduces us to a monomial situation, similar to that in
Example~\ref{SNC}. On deduces that $g \in \JJ(X, \lambda\cdot Y)$ if
and only if
$$\ord_{E_i} g \ge \lfloor\lambda a_i\rfloor - k_{i}$$
for every $i$. Equivalently, if we put
$\lfloor\lambda\mu^{-1}(Y)\rfloor=\sum_i\lfloor\lambda a_i \rfloor
E_i$, then
\begin{equation}\label{formula_multiplier}
\JJ(X, \lambda \cdot Y) = \mu_* \OO_{X'}(K_{X'/X} - \lfloor\lambda
\mu^{-1}(Y)\rfloor).
\end{equation}
 We refer to \cite{positivity} for details.

Note that because of the original definition, it follows that this
expression for $\JJ(X, \lambda \cdot Y)$ is independent of the
choice of a resolution of singularities. On the other hand, the
formula (\ref{formula_multiplier}) applies also when $X$ is non
necessarily affine. Note also that this formula implies that if
$\lambda_1\geq\lambda_2$, then
$$\JJ(X, \lambda_1\cdot Y) \subseteq \JJ(X, \lambda_2\cdot Y).$$

If $\lambda$ is small enough, then $\lambda a_i < k_i+1$ for $i =1,
\dots, N$. This implies that
$$ \ord_{E_i} 1 \ge \lfloor\lambda a_i\rfloor - k_{i},$$
hence $\JJ(X, \lambda\cdot Y) = \OO_X$. This leads us to the
definition of the log canonical threshold of the pair $(X, Y)$: this
is the smallest $\lambda$ such that $\JJ(X,\lambda\cdot Y)\neq
\OO_X$, i.e.
$$ c = \lc(X, Y) = \min_{i} \left\{\frac{k_i+1}{a_i}\right\}.$$
We may regard $\frac{1}{c}$ as a refined version of multiplicity. In
general a singularity with a smaller log canonical threshold tends
to be more complex.

The first appearance of the log canonical threshold was in the work
of Arnold, Gusein-Zade and Varchenko (see \cite{Arnold} and
\cite{Varchenko}), in connection with the behavior of certain
integrals over vanishing cycles. In the last decade this invariant
has enjoyed renewed interest due to its applications to birational
geometry. The following is probably the most interesting open
problem about log canonical thresholds.

\begin{conjecture}{\rm (}Shokurov{\rm )}
For every $n$, the set
$$\{\lc(X,Y)\mid \dim(X)=n, Y\subset X\}$$
satisfies the Ascending Chain Condition: it contains no strictly
increasing sequences.
\end{conjecture}

\smallskip

We can consider also higher jumping numbers. In general, we say that
$\lambda$ is a jumping number of $(X, Y)$, if
$$ \JJ(X, \lambda\cdot Y) \subsetneq \JJ(X, (\lambda - \epsilon)\cdot Y)$$
for all $\epsilon >0$. If $\lambda a_i$ is not an integer, then
$\lfloor\lambda a_i\rfloor=\lfloor(\lambda - \epsilon) a_i\rfloor$
for sufficiently small positive $\epsilon$. We see that a necessary
condition for $\lambda$ to be a jumping number is that $\lambda a_i$
is an integer for some $i$. In particular, if $\lambda$ is a jumping
number, then it is rational and has a bounded denominator.

The following theorem gives a periodicity property of the jumping
numbers.
\begin{theorem}
  \begin{enumerate}
  \item[(i)] If $Y = D$ is a hypersurface in $X$, then
  $$\JJ(X, \lambda\cdot
  D)\cdot \OO_X(-D) = \JJ(X, (\lambda+1)\cdot D).$$

  \item[(ii)] {\rm (Ein and Lazarsfeld \cite{EL2})} For every $Y$
  defined by the ideal $I_Y$, if
  $\lambda\ge\dim X-1$, then
    $$\JJ(X, \lambda\cdot Y)\cdot I_Y =\JJ(X, (\lambda+1)\cdot Y).$$
\end{enumerate}
\end{theorem}

\begin{corollary}
  If $\lambda > \dim X -1$, then $\lambda$ is a jumping number
  for $(X, Y)$ if and only if  so is $(\lambda+1)$.
\end{corollary}

We conclude that the set of jumping numbers of the pair $(X, Y)$ is
a discrete subset of $\QQ$ and it is eventually periodic with period
one.

\begin{example}
If $Y$ is a smooth subvariety of $X$ of codimension $e$, then the
set of jumping numbers of the pair $(X,Y)$ is $\{e, e+1, \cdots \}$.
 In particular $\lc(X,Y) = e$.
\end{example}

\begin{example}(Howald)
Let $X= \CC^n$ and let $Y$ be the closed subscheme defined by a
monomial ideal $\fra$. If $a = (a_1, a_2, \dots, a_n)\in\NN^n$, we
denote the monomial $x_1^{a_1}\cdots x_n^{a_n}$ by $x^a$. Consider
the Newton polyhedron $P_{\fra}$ associated with $\fra$: this is the
convex hull of those $a\in\NN^n$ such that $x^a\in\fra$. Using toric
geometry Howald showed in \cite{Howald} that
$$\JJ(X, \lambda\cdot Y)=(x^a\mid a+e\in {\rm Int}(\lambda P_{\fra})),$$
where $e=(1,\ldots,1)$. In particular, the log canonical threshold
$c$ of $(X,Y)$ is characterized by the fact $\frac{1}{c}\cdot e$
lies on the boundary of $P_{\fra}$.

For example, suppose that $\fra$ is the ideal $(x_1^{a_1},\ldots,
x_n^{a_n})$. In this case, the boundary of $P_{\fra}$ is
$$\{u=(u_1,\ldots,u_n)\in\RR_+^n\mid \sum_{i=1}^n\frac{u_i}{a_i}=
1\}.$$ Therefore $\lc(X, Y) = \sum_i\frac{1}{a_i}$.
\end{example}

\begin{example}
Suppose that $X = \CC^2$ and $Y$ is the plane cuspidal curve defined
by $x^3+y^5=0$. Then the set of jumping numbers for $Y$ is periodic
with period 1. The jumping numbers in $(0,1]$ are $\{\frac{8}{15},
\frac{11}{15}, \frac{13}{15}, \frac{14}{15},1\}$.
\end{example}

One reason that multiplier ideals have been very powerful in
studying questions in higher dimensional algebraic geometry is that
they appear naturally in a Kodaira type vanishing theorem. The
following statement is the algebraic version of a result due to
Nadel. In our context, it can be deduced from the Kawamata-Viehweg
Vanishing Theorem (see \cite{positivity}).

\begin{theorem} Let $X$ be a
  smooth projective variety and $Y$ a closed subscheme of
  $X$ defined by the ideal ${\mathcal I}_Y$.
If $A$ is a line bundle such that $I_Y\otimes A$ is globally
generated, and if $L$ is a line bundle such that $L-A$ is big and
nef, then for every $i>0$
  $$H^i(X, \OO_X(K_X+L) \otimes \JJ(X, Y)) = 0.$$
\end{theorem}

\section{Applications of multiplier ideals}

One of the most important applications of multiplier ideals is the
following theorem of Siu (see \cite{Siu1} and \cite{Siu2}) on the
deformation invariance of plurigenera.

\begin{theorem}
Let $f: X \lra T$ be a
  smooth projective morphism of relative dimension $n$
  between two smooth irreducible varieties. If we denote by $X_t$ the fiber
  $f^{-1}(t)$ for each $t \in T$, then for every fixed $m >0$,
  the dimension of the cohomology
  group $H^{0}(X_t, (\Omega_{X_t}^n)^{\otimes m})$ is
  independent of $t$.
\end{theorem}

The techniques involved in the proof of this theorem have been
recently applied by Siu, Hacon and McKernan to study one of the
outstanding problems in higher dimensional algebraic geometry, the
finite generation of the canonical ring (see, for example,
\cite{hacon}).

In a different direction, there are applications of multiplier
ideals to singularities of theta divisors on abelian varieties. Let
$(X, \Theta)$ be a principally polarized abelian variety, i.e.
$\Theta$ is an ample divisor on an abelian variety $X$ such that
$\dim H^0(X,\OO_X(\Theta)) = 1$. The following result is due to Ein
and Lazarsfeld \cite{EL}.

\begin{theorem}\label{PPAV}
Let $(X, \Theta)$ be a
principally polarized abelian
  variety. If $\Theta$ is irreducible, then $\Theta$ has at most
  rational singularities.
\end{theorem}

\begin{corollary} Let~$(X,\Theta)$ be a principally polarized abelian
variety of dimension $g$, with $\Theta$ irreducible. If
$$\Sigma_k(\Theta) = \{x \in X\mid  \mult_x(\Theta)\ge k\},$$
then for every $k \ge 2$ we have $\codim(\Sigma_k(\Theta),X) \ge
k+1$. In particular, $\Theta$ is a normal variety and
$\mult_x(\Theta) \le g-1$ for every singular point $x$ on $\Theta$.
\end{corollary}

\begin{remark}
  The fact that $\Theta$ is normal was first conjectured by Arbarello,
  De Concini and Beauville. When $X$ is the Jacobian of a curve, the
  fact that $\Theta$ has only rational singularities was proved by
  Kempf. Note also that in this case, a classical theorem of Riemann expresses
  the multiplicity of $\Theta$ at a point in term of the dimension of
  the
  corresponding linear system on the curve. It was  Koll\'{a}r
  who first observed in \cite{Sh}
  that one can use vanishing theorems to study the singularities of the theta
  divisor: he showed that for every principally polarized abelian
  variety $(X, \Theta)$, we have $\lc(X,\Theta)=1$.
  Theorem~\ref{PPAV} above is a strengthening of Koll\'{a}r's
  result.
\end{remark}

Multiplier ideals have  been applied in several other directions: to
Fujita's problem on adjoint linear systems \cite{AS}, to Effective
Nullstellensatz \cite{EL2}, to Effective Artin-Rees Theorem
\cite{ELSV}. Building on work of Tsuji, recently Hacon and McKernan
and independently, Takayama have used multiplier ideals to prove a
very interesting result on boundedness of pluricanonical maps for
varieties of general type (see \cite{hacon1} and \cite{Takayama}).
We end this section with an application to commutative algebra due
to Ein, Lazarsfeld and Smith \cite{ELS}.

Let $X$ be a smooth $n$-dimensional variety and $Y\subseteq X$
defined by the reduced sheaf of ideals $\fra$. The $m^{\rm th}$
symbolic power of $\fra$ is the sheaf $\fra^{(m)}$ of functions on
$X$ that vanish with multiplicity at least $m$ at the generic point
of every irreducible component of $Y$. If $Y$ is smooth, then the
symbolic powers of $\fra$ agree with the usual powers, but in
general they are very different.

\begin{theorem}
If $X$ is a smooth $n$-dimensional variety
and if $\fra$ is a reduced sheaf of ideals, then
$\fra^{(mn)}\subseteq
\fra^m$ for every $m$.
\end{theorem}

\vfill\eject

\section{Bounds on log canonical thresholds and birational rigidity}

In this section we compare the log canonical threshold with the
classical Samuel multiplicity. We give then an application of the
inequality between these two invariants to a classical question on
birational rigidity. Let $X$ be a smooth complex variety and $x \in
X$ a point. Denote by $R$ the local ring of $X$ at $x$, and by
${\frak{m}}$ its maximal ideal. The following result was proved by
de Fernex, Ein and Musta\c{t}\v{a} in \cite{DEM}.

\begin{theorem}\label{ineq_lc}
Let $\fra$ be an ideal in $R$ that defines a subscheme $Y$ supported
at $x$. Let $c$ be the log canonical threshold of $(X,Y)$,
$l(R/\fra)$ be the  length of $R/\fra$ and $e(\fra)$ be the Samuel
multiplicity of $R$ along $\fra$. If $ n = \dim R$, then we have the
following inequalities.
\begin{enumerate}
\item[(i)] $l(R/\fra) \ge \frac{n^{n}}{n!\cdot c^n}$.
\item[(ii)] $e(\fra) \ge \frac{n^{n}}{c^n}.$ Furthermore,
this is an equality if and only if the integral closure of $\fra$ is
equal to $\frak{m}^k$ for some $k$.
\end{enumerate}
\end{theorem}

The first assertion in (ii) above can be easily deduced from (i).
The proof of (i) proceeds by reduction to the monomial case, via a
Gr\"{o}bner deformation. When $\fra$ is monomial,  the inequality
follows by a combinatorial argument from the explicit description of
the invariants.

\begin{example}
 Suppose that $\fra = (x_1^{a_1}, \dots, x_n^{a_n})$.
 In this case $e(\fra) = \prod_{i=1}^n a_i$ and $\lc(\fra) =
 \sum_{i=1}^n \frac{1}{a_i}$. The inequality in
 Theorem~\ref{ineq_lc}(ii) becomes
 $$ \prod_{i=1}^n a_i \ge \frac{n^{n}}{\left(\sum_{i=1}^n \frac{1}{a_i}\right)^n}.$$
This is equivalent to
$$ \left(\frac{1}{n} \sum_{i=1}^n \frac{1}{a_i}\right)^n \ge
\prod_{i=1}^n \frac{1}{a_i},$$ which is just the classical
inequality between the arithmetic and the geometric mean.
\end{example}

\begin{remark} When $X$ is a surface,
the inequality in (ii) above was first proved by Corti \cite{Corti}.
\end{remark}

Theorem~\ref{ineq_lc} is used in \cite{DEM2} to study the behavior
of the log canonical threshold under a generic projection. More
generally, one proves the following

\begin{theorem}
Let $f: X \lra Y$ be a smooth proper morphism of relative dimension
$k-1$ between two smooth complex varieties. If $V$ is a locally
complete intersection of codimension $k$ in $X$ such that $f|_{V}$
is finite, then
$$\lc(Y, f_*(V)) \le\frac{\lc(X,V)^k}{k^{k}}.$$
\end{theorem}

Using the above theorems and some beautiful geometric ideas of
Pukhlikov \cite{Pu1}, one gives in \cite{DEM2} a simple uniform
proof for the following result.

\begin{theorem}\label{thm10}
If $X$ is a smooth hypersurface of degree $N$ in $\CC\PP^N$, with $4
\le N \le 12$, then $X$ is birationally superrigid. In particular,
every birational automorphism of $X$ is biregular.
\end{theorem}

\begin{remark}
Consider the group $\Aut_{\CC}(\CC(X))$, the automorphism group of
the field of rational functions of $X$. This is naturally isomorphic
to $\Bir_{\CC}(X)$, the group of birational automorphisms of $X$. If
$X$ is birationally superrigid, then $\Bir_{\CC}(X) \simeq
\Aut_{\CC}(X)$, the automorphism group of $X$. When $X$ is a
hypersurface of degree $N$ in $\PP^N$,  $X$ has no nonzero vector
fields and therefore $\Aut_{\CC}(X)$ is a finite group. This shows
that $X$ is not a rational variety: if $\CC(X)$ is purely
transcendental, then $\Aut_{\CC}(X)$ will contain a subgroup
isomorphic to the general linear group $GL_n$.
\end{remark}

\begin{remark}
When N=4, it is a classical theorem of Iskovskikh and Manin that $X$
is birationally rigid \cite{IM}. They used this to show that the
function field of a suitable quartic threefold provides a
counterexample to the classical Luroth's problem: $\CC(X)$ is a
$\CC$-subfield of the purely transcendental field $\CC(x_1, x_2,
x_3)$ but $\CC(X)$ is not purely transcendental. For N=5,
Theorem~\ref{thm10} is a result of Pukhlikov \cite{Pu2}. The cases
$N= 6$, $7$ and $8$ were first proved by Cheltsov \cite{Cheltsov}.
We mention also that Pukhlikov \cite{Pu5} has shown that a generic
hypersurface of degree $N$ in $\PP^N$ is superrigid for every $N \ge
4$.
\end{remark}

\section{Bernstein-Sato polynomials}

Let $f \in \CC[x_1, x_2, \ldots, x_n]$ be a nonzero polynomial. We
denote by $A_n$ the Weyl algebra of differential operators on
${\mathbb A}^n$, that is
$$A_n=\CC[x_1,\ldots,x_n,\partial_{x_1},\ldots,\partial_{x_n}].$$ Let $s$ be another
variable and consider the following functional equation:
\begin{equation}\label{bernstein1}
b(s)f^s=P(s,x,\partial_{x})\bullet f^{s+1},
\end{equation}
 where $b(s) \in \CC[s]$ and $P\in A_n[s]$. Here $f^s$ is considered
 a formal symbol, and the action $\bullet$ of $P$ is defined
 via $\partial_{x_i}\bullet f^s=sf^{-1}\frac{\partial
 f}{\partial x_i}f^s$. On the other hand, if we let $s=m$ for an
 integer $m$, then (\ref{bernstein1}) has the obvious meaning.

It is easy to see that the set of polynomials $b(s)$ for which there
is $P$ satisfying (\ref{bernstein1}) is an ideal in the polynomial
ring $\CC[s]$. It was proved by Bernstein in \cite{Bernstein} using
the theory of holonomic $D$-modules that this ideal is nonzero. Its
monic generator is denoted by $b_f(s)$ and is called the
Bernstein-Sato polynomial of $f$. It is an interesting and subtle
invariant of the singularities of the hypersurface defined by $f$.

\begin{example}
\begin{enumerate}
\item Making $s=-1$ in (\ref{bernstein1}) we see that $b_f(-1)f^{-1}$
lies in $\CC[x_1,\ldots,x_n]$. If $f$ is nonconstant, it follows
that $-1$ is a root of $b_f$.

\item If $f=x$, then $b_f(s)=(s+1)$. Indeed, we have
$$(s+1)f^s=\partial_{x}\bullet f^{s+1}.$$
More generally, if $f$ defines a nonsingular hypersurface, then
$b_f(s)=(s+1)$.

\item If $f=x_1^2+\ldots+x_n^2$, then
$b_f(s)=(s+1)\left(s+\frac{n}{2}\right)$ and
$$b_f(s)f^s=\frac{1}{4}(\partial_{x_1}^2+\ldots+\partial_{x_n}^2)\bullet
f^{s+1}.$$

\item If $f=x^2+y^3$, then
$b_f(s)=\left(s+\frac{5}{6}\right)(s+1)\left(s+\frac{7}{6}\right)$
and
$$b_f(s)f^s=\left(\frac{1}{27}\partial_y^3+\frac{1}{6}y\partial_x^2\partial_y+\frac{1}{8}x\partial_x^3+\frac{3}{8}\partial_x^2\right)
\bullet f^{s+1}.$$
\end{enumerate}
\end{example}

Computing Bernstein-Sato polynomials in general is quite subtle (see
\cite{Yano}). On the other hand, there has been a lot of recent
progress in algorithmic computation using Gr\"{o}bner bases in the
Weyl algebra (see \cite{SST}).

We describe now the connection between the roots of the
Bernstein-Sato polynomial of $f$ and the jumping numbers of the
hypersurface $Y$ defined in $\AAA^n$ by $f$. An important theorem of
Kashiwara \cite{Kashiwara} asserts that all the roots of $b_f(s)$
are negative rational numbers. Building on Kashiwara's work, Lichtin
made this more explicit in \cite{Lichtin}, describing a connection
between the roots of $b_f(s)$ and a log resolution of the pair
$(\AAA^n,Y)$. This says that if $\mu\colon X'\to\AAA^n$ is a log
resolution of $(X,Y)$, then every root of $b_f(s)$ is of the form
$-\frac{k_i+m}{a_i}$ for some $i$ and some positive integer $m$ (we
use the notation introduced in \S 2). In particular, we see that
every root of $b_f(s)$ is rational, and no larger than
$-\lc(\AAA^n,Y)$. However, we stress that unlike in the case of
multiplier ideals, there is no explicit description of the
Bernstein-Sato polynomial in terms of a log resolution.

On the other hand, the following result of Ein, Lazarsfeld, Smith
and Varolin \cite{ELSV} shows that in a suitable range, all jumping
numbers give roots of the Bernstein-Sato polynomial.

\begin{theorem}\label{bernstein2}
  If $\lambda\in (0,1]$ is a jumping number of $(\AAA^n,Y)$, then
  $-\lambda$ is a root of the Bernstein-Sato polynomial $b_f(s)$.
\end{theorem}

The proof of this theorem uses the functional equation
(\ref{bernstein1}) and integration by parts. The case when
$\lambda=\lc(\AAA^n,Y)$ was proved by Koll\'{a}r in \cite{Kollar1}.
Note that in conjunction with the above mentioned result of Lichtin,
this gives the following

\begin{corollary}
The largest root of $b_f(s)$ is $-\lc(\AAA^n,Y)$.
\end{corollary}

\bigskip

A different point of view on the connection between multiplier
ideals and Bernstein-Sato polynomials was given by Budur and Saito.
In fact, they show how to recover the multiplier ideals from a
filtration that appears in $D$-module theory, the $V$-filtration. We
present now their result.

Let $t$ be a new variable, and let $A_{n+1}$ denote the Weyl algebra
corresponding to the affine space $\AAA^{n+1}$, with coordinates
$x_1,\ldots,x_n,t$. We consider the module $B_f$ that is the first
local cohomology module of $\AAA^{n+1}$ along the embedding of
$\AAA^n$ as the graph of $f$, i.e.
$$B_f=\CC[x_1,\ldots,x_n,t]_{f-t}/\CC[x_1,\ldots,x_n,t].$$
Let $\delta$ be the class of $\frac{1}{f-t}$ in $B_f$ ($\delta$ is
the "delta-function associated to the graph of $f$").

Note that $B_f$ has a natural structure of left module over
$A_{n+1}$. Since $\partial_t^m\delta$ is the class of
$\frac{m!}{(f-t)^{m+1}}$ in $B_f$, we see that $B_f$ is free over
$\CC[x_1,\ldots,x_n]$, with basis given by $\{\partial_t^j
\delta\mid j\geq 0\}$.

The $V$-filtration is a decreasing filtration on $B_f$ by left
$A_n$-submodules $V^{\alpha}$ indexed by $\alpha\in\QQ$, with the
following properties:
\begin{enumerate}
\item[(i)] $\bigcup_{\alpha}V^{\alpha}=B_f$.

\item[(ii)] The filtration is semicontinuous and discrete in the following
sense: there is a positive integer $\ell$ such that for every
integer $m$ and every $\alpha\in
\left(\frac{m-1}{\ell},\frac{m}{\ell}\right]$ we have
$V^{\alpha}=V^{m/\ell}$.

\item[(iii)]  We have $t\cdot V^{\alpha}\subseteq V^{\alpha+1}$ for every
$\alpha$, with equality if $\alpha>0$.

\item[(iv)] We have $\partial_t\cdot V^{\alpha}\subseteq V^{\alpha-1}$
for every $\alpha$.

\item[(v)] For every $\alpha$, if we put
$V^{>\alpha}:=\cup_{\beta>\alpha}V^{\beta}$, then
$(\partial_tt-\alpha)$ is nilpotent on $V^{\alpha}/V^{>\alpha}$.
\end{enumerate}

The key property is (v) above. One can think of the $V$-filtration
as an attempt to diagonalize the operator $\partial_tt$ on $B_f$. It
is not hard to show that if a filtration as above exists, then it is
unique. Malgrange \cite{Malgrange} proved the existence of the
$V$-filtration using the existence of the Bernstein-Sato polynomial
and the rationality of its roots. To explain the role played by
$b_f(s)$ in the construction of the $V$-filtration, we mention that
the equation (\ref{bernstein1}) in the definition of $b_f$ is
equivalent with the following equality in $B_f$:
\begin{equation}
b(-\partial_tt)\cdot\delta=P(-\partial_tt,x,\partial_x)\cdot
t\delta.
\end{equation}

The following result of Budur and Saito \cite{BS} shows that the
multiplier ideals can be obtained as a piece of the $V$-filtration.
We consider $\CC[x_1,\ldots,x_n]$ embedded in $B_f$ by $h\to
h\delta$.

\begin{theorem}\label{bernstein4}
If $Y$ is the hypersurface defined by $f$, then for every
$\lambda>0$ we have $\JJ(\AAA^n,\lambda\cdot
Y)=V^{\lambda+\epsilon}\cap \CC[x_1,\ldots,x_n]$, where
$0<\epsilon\ll 1$.
\end{theorem}

The assertion in Theorem~\ref{bernstein2} can be deduced from this
statement. The proof of Theorem~\ref{bernstein4} involves two steps.
First, one describes the $V$-filtration in the case when $f$ defines
a divisor with simple normal crossings: $f=x_1^{a_1}\ldots
x_n^{a_n}$. In this case, let us put $\JJ'(\AAA^n,\alpha\cdot
Y):=\JJ(\AAA^n,(\alpha-\epsilon)\cdot Y)$ for $0<\epsilon\ll 1$
(with the convention $\JJ'(\AAA^n,\alpha\cdot
Y)=\CC[x_1,\ldots,x_n]$ if $\alpha\leq 0$). If we take
 $V^{\alpha}$ to be generated over $A_n$ by $\JJ'(\AAA^n,(\alpha+j)\cdot Y)\partial_t^j\delta$, where $j$ varies
over the nonnegative integers, then one can check that these
$V^{\alpha}$ satisfy the properties in the definition of the
$V$-filtration. In particular, this easily implies the statement of
Theorem~\ref{bernstein4} in this case. The hard part of the proof
uses Saito's theory of mixed Hodge modules to deduce the general
case of the theorem by relating the $V$-filtrations of $f$ and of a
log resolution.

\smallskip

We mention that Kashiwara constructed in \cite{Kashiwara2} a
$V$-filtration associated to several polynomials.  Budur,
Musta\c{t}\v{a} and Saito used this in \cite{BMS} to introduce and
study the Bernstein-Sato polynomial associated to a subscheme not
necessarily of codimension one, and to generalize
Theorems~\ref{bernstein2} and \ref{bernstein4} to this setting.

\vfill\eject

\section{Spaces of arcs and contact loci}

Let $X$ be a smooth $n$-dimensional complex variety. Given $m \ge
0$, we denote by
$$  X_m \ = \ \text{Hom} \big( \,  \Spec \, \CC[t]/
(t^{m+1}) \, , \, X \, \big) $$ the space of $m^{\text{th}}$ order
jets on $X$. This carries a natural scheme structure. Similarly we
define the space of formal arcs on $X$ as
$$ X_{\infty} = \text{Hom} \big(\, \Spec \, \CC[[t]]
\, , X \, \big).$$ These constructions are functorial, hence to
every morphism $\mu\colon X'\to X$ we associate corresponding
morphisms $\mu_m$ and $\mu_{\infty}$. Thanks to the work of
Kontsevich, Denef, Loeser and others on motivic integration, in
recent years these spaces have been very useful in constructing
invariants of singular algebraic varieties. In what follows we
describe some applications of these spaces to the study of
singularities.

We have natural projection maps induced by truncation $X_{m+1}\to
X_m$. Since $X$ is smooth, this is locally trivial in the Zariski
topology, with fiber $\AAA^n$. We similarly have projection maps
$X_{\infty}\to X_m$. A subset $C$ of $X_{\infty}$ is called a
\emph{cylinder} if it is the inverse image of a constructible set
$S$ in some $X_m$. Moreover, $C$ is called locally closed (closed,
irreducible) if $S$ is so. If $C$ is a closed cylinder that is the
inverse image of $S\subset X_m$, its codimension in $X_{\infty}$ is
equal to the codimension of $S$ in $X_m$.

Consider a nonzero ideal sheaf $\fra \subseteq \OO_X$ defining a
subscheme $Y \subset X$.  Given a finite jet or an arc $\gamma$ on
$X$, the
 \text{order of vanishing}
 of $\fra$  --- or
the \text{order of contact} of the corresponding
scheme $Y$ --- along
$\gamma$ is defined in the natural
way. Specifically,  pulling $\fra $ back  via
$\gamma$ yields an ideal $(t^e)$ in $\CC[t]/(t^{m+1})$ or
$\CC[[t]]$, and one sets
\[  \ord_\gamma(\fra) \ = \ \ord_\gamma(Y) \ = \ e.
\] (Take  $\ord_\gamma(\fra) = m+1$ when
$\fra$ pulls back to the zero ideal in $\CC[t]/(t^{m+1})$ and
$\ord_\gamma(\fra) = \infty$ when it pulls back to the zero ideal in
$\CC[[t]]$.) For a fixed integer $p \ge 0$, we define the
\textit{contact locus}
$$ \Cont^p(Y) \ = \ \Cont^p(\fra) \ =_{\text{def}} \ \big\{ \, \gamma
\in X_{\infty} \, | \, \ord_{\gamma}(\fra) = p \, \big \}. $$ Note
that this is a locally closed cylinder: for $m\geq p$, it is the
inverse image of
\begin{equation}
 \Cont^p(Y)_m \ = \  \Cont^p(\fra)_m \ =_{\text{def}} \
\big \{ \,\gamma \in X_m \mid \ \ord_\gamma(\fra) = p \, \big \},
\end{equation}
which is locally closed in $X_m$. A subset of $X_{\infty}$ is called
an \emph{irreducible closed contact subvariety} if it is the closure
of an irreducible component of $\Cont^p(Y)$ for some $p$ and $Y$.

Suppose now that $W$ is an arbitrary irreducible closed cylinder in
$X_{\infty}$. We can naturally associate a valuation of the function
field of $X$ to $W$ as follows. If $f$ is a nonzero rational
function of $X$, we put
$$ \val_W(f) = \ord_\gamma(f)  \ \ \text{for a general $\gamma \in W$}. $$
This valuation is not identically zero if and only if $W$ does not
dominate $X$.

If $\mu : X^\pr \lra X$ is a proper birational morphism, with
$X^\pr$ smooth, and if $E$ is an irreducible divisor on $X^\pr$,
then we define a valuation by
$$\val_E(f) = \ \ \text{the vanishing order of $f\circ\mu$ along $E$}.$$
A valuation on the function field of $X$ is called a
\emph{divisorial valuation} (with center on $X$) if it is of the
form $m\cdot \val_E$ for some positive integer $m$ and some divisor
$E$ as above.

A key invariant associated to a divisorial valuation $v$ is its
\emph{log discrepancy}. If $E$ is a divisor as above, we put $k_E=
\val_E({\rm det}(J(\mu))$, where
 $J(\mu)$ is the Jacobian matrix of $\mu$. Equivalently, $k_E$
 is the coefficient of $E$ in the relative canonical divisor
 $K_{X'/X}$. Note that $k_E$
depends only on ${\rm val}_E$ (it does not depend on the model
$X'$). Given an arbitrary divisorial valuation
 $m\cdot\val_E$, we define its log discrepancy as $m(k_E+1)$.

Consider a divisor $E$ on $X'$ as above. If $C_m(E)$ is the closure
of $\mu_{\infty}({\rm Cont}^m(E))$, then it is not hard to see that
$C_m(E)$ is an irreducible closed contact subvariety of $X_{\infty}$
such that $\val_{C_m(E)}=m\cdot\val_E$. The following result of Ein,
Lazarsfeld and Musta\c{t}\v{a} \cite{ELM} describes in general the
connection between cylinders and divisorial valuations.

\begin{theorem}\label{valuations}
Let $X$ be a smooth variety.
\begin{enumerate}
\item[(i)] If $W$ is an irreducible, closed cylinder in $X_{\infty}$
that does not dominate $X$, then the valuation $\val_W$ is
divisorial.
\item[(ii)] For every divisorial valuation $m\cdot\val_E$, there is
a unique maximal irreducible closed cylinder $W$ such that
$\val_W=m\cdot\val_E$: this is $W=C_m(E)$.
\item[(iii)]
The map that sends $m\cdot \val_E$ to $C_m(E)$ gives a bijection
between divisorial valuations of $\CC(X)$ with center on $X$ and the
set of irreducible closed contact subvarieties of $X_{\infty}$.
\end{enumerate}
\end{theorem}

The applicability of this result to the study of singularities is
due to the following description of log discrepancy of a divisorial
valuation in terms of the codimension of a certain set of arcs.

\begin{theorem}\label{geometric}
Given a divisorial valuation $v=m\cdot \val_E$ with center on $X$,
if $C_m(E)$ is its associated irreducible closed contact subvariety
in $X_{\infty}$, then the log discrepancy of $v$ is equal to
$\codim(C_m(E), X_{\infty})$.
\end{theorem}

Combining the statements of the above theorems, we deduce a lower
bound for the codimension of an arbitrary cylinder in terms of the
log discrepancy of the corresponding divisor.

\begin{corollary}
If $W$ is a closed, irreducible cylinder in $X_{\infty}$ that does
not dominate $X$, then $\codim(W,X_{\infty})$ is bounded below by
the log discrepancy of $\val_W$.
\end{corollary}

\begin{remark}
The above two theorems also hold for singular varieties after some
minor modifications using Nash's blow-up and Mather's canonical
class.
\end{remark}

The key ingredient in the proof of the above theorems is the
following result due to Kontsevich, Denef and Loeser (see
\cite{denef}. It constitutes the geometric content of the so-called
Change of Variable Theorem in motivic integration. Suppose that $\mu
: X^\pr \lra X$ is a proper, birational morphism of smooth varieties
and let $K_{X^\pr / X}$ be the relative canonical divisor.

\vfill\eject

\begin{theorem}
Given integers $e \ge 0$ and $m\geq e$, consider the contact locus
\[ \Cont^e(K_{X^\pr/X})_m \ = \ \big\{\, \gamma^\pr \in
X^\pr_m  \mid \ord_{\gamma^\pr}(K_{X^\pr/X}) = e \, \big \}. \] If
$m \ge 2e$, then $\Cont^e(K_{X^\pr/X})_m $ is a union of fibres of
$\mu_m : X_m^\pr \lra X_m$, each of which is isomorphic to an affine
space $\AAA^e$. Moreover, if \[ \gamma^\pr \ , \ \gamma^{\pr \pr} \
\in \ \Cont^e(K_{X^\pr/X})_m \]
 lie in the same fibre of
$\mu_m$, then they have the same image in
$X^\pr_{m-e}$.
\end{theorem}

As an  application of Theorems~\ref{valuations} and \ref{geometric},
one gives in \cite{ELM} a simple proof of the following result of
Musta\c{t}\v{a} \cite{Mustata1} describing the log canonical
threshold in terms of the geometry of the space of jets.

\begin{theorem}
Let $X$ be a smooth complex variety and $Y$ be a closed subscheme of
$X$ defined by the nonzero ideal sheaf ${\mathcal I}_Y$. Let $X_m$
and $Y_m$ be the spaces of $m^{\rm th}$ order jets of $X$ and $Y$,
respectively. If $c = \lc(X, Y)$, then
  \begin{enumerate}
    \item[(i)] For every $m$ we have $\codim(Y_m, X_m) \ge c\cdot (m+1)$.
    More generally, if $W\subset X_{\infty}$ is an irreducible closed
    cylinder that does not dominate $X$, then
    ${\rm codim}(W,X_{\infty})\geq c\cdot{\rm val}_W({\mathcal I}_Y)$.
    \item[(ii)] If $m$ is sufficiently divisible, then
      $\codim(Y_m, X_m) = c\cdot (m+1)$.
    \item[(iii)] We have $ c = \lim_{m \rightarrow \infty} \frac{\codim(Y_m, X_m)}{m+1}.$
   \end{enumerate}
\end{theorem}

The above results relating divisorial valuations with the space of
arcs can be used to study more subtle invariants of singularities of
pairs. Let $Y$ be a closed subscheme of the smooth variety $X$, and
let $\lambda$ be a positive real number. We associate a numerical
invariant to the pair $(X,\lambda\cdot Y)$ and to an arbitrary
nonempty closed subset $B\subseteq X$, as follows.

Consider a divisorial valuation of the form $\val_E$ with center
$c_X(E)$ in $X$ (the center is the image of $E$ in $X$).  The
\emph{log discrepancy} of the pair $(X, \lambda\cdot Y)$ along $E$
is
\[ a(E, X, \lambda\cdot Y) = k_E + 1 - \lambda \cdot \val_E(I_Y), \]
where $I_Y$ is the ideal of $Y$ in $X$. Note that if $Y=\emptyset$,
we recover the log discrepancy of $\val_E$. The idea is to measure
the singularities of the pair $(X, \lambda\cdot Y)$  using the log
discrepancies along divisors with center contained in $B$.

\begin{definition}\label{def1}
Let $B\subset X$ be a nonempty closed subset. The minimal log
discrepancy of $(X,\lambda \cdot Y)$ over $B$ is defined by
\begin{equation}
{\rm mld}(B;X,\lambda \cdot Y):=\inf_{c_X(E)\subseteq B}\{a(E; X,
\lambda \cdot Y)\}.
\end{equation}
\end{definition}

\begin{remark}
One can show that ${\rm mld}(B;X,\lambda \cdot Y)$ is either
$-\infty$ or a nonnegative real number. In fact, ${\rm
mld}(B;X,\lambda\cdot Y)\neq-\infty$ if and only if there is an open
neighborhood $U$ of $B$ such that $\lc(U,U\cap Y)\geq\lambda$. An
important fact about minimal log discrepancies is that they can be
computed using a log resolution of $(X,B\cup Y)$, see \cite{Ambro}.
\end{remark}

The following theorem of Ein, Musta\c{t}\v{a} and Yasuda \cite{EMY}
gives a description of minimal log-discrepancies in term of the
geometry of the space of arcs.

\begin{theorem} Let $B$ be a nonempty, proper closed subset of $X$,
and let
 $\pi\colon X_{\infty} \lra X$ be the projection map.
 For every proper closed subscheme $Y$ of $X$ and for every $\lambda$
 and $\tau\in\RR_+$ we have
${\rm mld}(B;X,\lambda \cdot Y) \ge \tau$ if and only if for every
irreducible closed cylinder $W\subseteq\pi^{-1}(B)$ we have
$$\codim(W, X_{\infty}) \ge \lambda \cdot \val_W(I_Y) +\tau.$$
\end{theorem}

The above theorem can be applied to study the behavior of
singularities of pairs under restriction to a divisor. This is
useful whenever one wants to do induction on dimension. Suppose that
$D$ is a smooth divisor on $X$. We want to relate the singularities
of $(X,\lambda\cdot Y)$ with those of $(D,\lambda\cdot Y\vert_D)$.
The adjunction formula suggests that the precise relation should be
between $(X,D+\lambda\cdot Y)$ and $(D,\lambda\cdot Y\vert_D)$. The
precise formula is the content of the following theorem from
\cite{EMY}.

\begin{theorem}\label{Inversion}
Let $D$ be a smooth divisor on the smooth variety $X$ and let $B$ be
a nonempty proper closed subset of $D$. If $Y$ is a closed subscheme
of $X$ such that $D\not\subseteq Y$, and if $\lambda\in\RR_+$, then
$${\rm mld}(B;X,D+\lambda \cdot Y)={\rm mld}(B;D,\lambda \cdot Y\vert_D).$$
\end{theorem}

\begin{remark}
The notion of minimal log discrepancy plays an important role in the
Minimal Model Program. It can be defined under weak assumptions on
the singularities of $X$: one requires only that $X$ is normal and
$\QQ$-Gorenstein. Koll\'{a}r and Shokurov have conjectured the
statement of Theorem~\ref{Inversion} with the assumption that $X$
and $D$ are only normal and $\QQ$-Gorenstein. It is easy to see that
the inequality "$\leq$" holds in general, and the opposite
inequality is known as Inversion of Adjunction (see \cite{Kollar1}
and \cite{Kollar2} for a discussion of this conjecture and related
topics). Theorem~\ref{Inversion} has been generalized in \cite{EM}
to the case when both $X$ and $D$ are normal locally complete
intersections.
\end{remark}

The interpretation of minimal log discrepancies in terms of spaces
of arcs gives also the following semicontinuity statement. This was
conjectured for an arbitrary (normal and $\QQ$-Gorenstein) variety
$X$ by Ambro and Shokurov, see \cite{Ambro}. The statement below,
due to Ein, Musta\c{t}\v{a} and Yasuda \cite{EMY} has been
generalized to the case of a normal locally complete intersection
variety in \cite{EM}.

 \begin{theorem}\label{semicont}
If $X$ is a smooth variety and if $Y$ is a closed subscheme of $X$,
then for every $\lambda\in\RR_+$, the function on $X$ defined by
$x\longrightarrow{\rm mld}(x;X,\lambda \cdot Y)$ is lower
semicontinuous.
\end{theorem}

We end with a result that translates properties of the minimal log
discrepancy over the singular locus of a locally complete
intersection variety into geometric properties of its spaces of
jets.

\begin{theorem}\label{singX_m}
Let $X$ be a normal locally complete intersection variety of
dimension $n$.
\begin{enumerate}
 \item[(i)] $X_m$ has pure
dimension $n(m+1)$ for every $m$ {\rm (}and in this case $X_m$ is
also a locally complete intersection{\rm )} if and only if ${\rm
mld}(X_{\rm sing};X,\emptyset)\geq 0$ {\rm (}this says that $X$ has
\emph{log canonical singularities}{\rm )}.
 \item[(ii)]  $X_m$ is irreducible for every $m$
 (and in this case it is also reduced) if and only if
 ${\rm mld}(X_{\rm sing};X,\emptyset)\geq 1$ {\rm (}this
 says that $X$ canonical singularities{\rm )}.
 \item[(iii)] $X_m$ is normal for every $m$ if and only if
 ${\rm mld}(X_{\rm sing};X,\emptyset)>1$ {\rm (}this says that $X$
 has \emph{terminal singularities}{\rm )}.
 \item[(iv)] In general, we have
  $\codim((X_m)_{\rm sing}, X_m) \ge {\rm mld}(X_{\rm sing};X,\emptyset)$ for every $m$.
\end{enumerate}
\end{theorem}

\begin{remark}  The description in (ii) above was first proved in
\cite{Mustata2}. Note that since $X$ is in particular Gorenstein, it
is known that $X$ has canonical singularities if and only if it has
rational singularities. All the statements in the above theorem were
obtained in \cite{EMY} and \cite{EM} combining the description of
minimal log discrepancies in terms of spaces of arcs and Inversion
of Adjunction.
\end{remark}

\section{Invariants in positive characteristic}

Several invariants have been recently introduced in positive
characteristic using the Frobenius morphism, invariants whose
behavior is formally very similar to the ones we have discussed in
characteristic zero. Moreover, there are interesting results and
conjectures involving the comparison between the two sets of
invariants via reduction mod $p$.

As in the case of singularities of pairs $(X,Y)$ in characteristic
zero, one can develop the theory under very mild assumptions on the
ambient variety $X$ (in fact, the positive characteristic theory
does not even need the assumption that $X$ is $\QQ$-Gorenstein). For
this one needs to use the full power of the theory of tight closure
of Hochster and Huneke \cite{HH}. However, the definitions become
particularly transparent if we assume $X$ nonsingular. Therefore, in
accord with the setup in the previous sections, we will make this
assumption. The theory we present here is due to Hara and Yoshida
\cite{HY} building on previous work of Hara, Smith, Takagi and
Watanabe.

We work in the local setting with a regular local ring
$(R,\frakm,k)$ of characteristic $p>0$. Let $n=\dim(R)$ and let $E$
be the top local cohomology module of $R$, $E=H_{\frakm}^n(R)$. If
$x_1,\ldots,x_n$ generate $\frakm$, then
\begin{equation}\label{loc_coh}
E\simeq R_{x_1\ldots
x_n}/\sum_{i=1}^nR_{x_1\ldots\widehat{x_i}\ldots x_n}.
\end{equation}
The Frobenius morphism on $R$ induces a Frobenius morphism $F_E$ on
$E$ that via the isomorphism (\ref{loc_coh}) takes the class of
$u/(x_1\cdots x_n)^d$ to the class of $u^p/(x_1\cdots x_n)^{pd}$.

We want to study the singularities of the pair $(X,Y)$, where
$X={\rm Spec}(R)$ and $Y$ is defined by a nonzero ideal $\fra$. For
every $r\geq 0$ and every $e\geq 1$, we put
$$Z_{r,e}:={\rm ker}(\fra^rF_E^e)=\{w\in E\mid h F_E^e(w)=0\,{\rm
for}\,{\rm all}\,h\in\fra^r\}.$$ Given a nonnegative real number
$\lambda$, the \emph{test ideal} of the pair $(X,\lambda\cdot Y)$ is
$$\tau(X,\lambda\cdot Y):={\rm Ann}_R\left(\bigcap_{e\geq
1}Z_{\lceil \lambda p^e\rceil,e}\right).$$ Here $\lceil\alpha\rceil$
denotes the smallest integer that is $\geq\alpha$.

As Hara and Yoshida show in \cite{HY}, the test ideals
$\tau(X,\lambda\cdot Y)$ enjoy formal properties similar to those of
the multiplier ideals $\JJ(X,\lambda\cdot Y)$ in characteristic
zero. In particular, we can consider the jumping numbers for the
test ideals: these are the $\lambda$ such that $\tau(X,\lambda\cdot
Y)\subsetneq\tau(X,(\lambda-\epsilon)\cdot Y)$ for every positive
$\epsilon$.

The set of jumping numbers for the test ideals are also eventually
periodic with period one. However, two basic properties that for
multiplier ideals follow simply from the description in terms of a
log resolution are not known for test ideals: it is not known
whether every jumping number for the test ideals is rational, and
whether in every bounded interval there are only finitely many such
jumping numbers. We want to stress that the problem does \emph{not}
come from the fact that we do not know, in general, whether such
resolutions exist. Even when we have such resolutions, the
invariants in characteristic $p$ do not depend simply on the
numerical data of the resolution (see Example~\ref{cusp} below for
the case of the cusp).

\smallskip

There is a more direct description of the set of jumping numbers
given by Musta\c{t}\v{a}, Takagi and Watanabe in \cite{MTW}. Suppose
that $J$ is a proper ideal of $R$ containing $\fra$ in its radical.
For every $e\geq 1$, define $\nu^J(p^e)$ to be the largest $r$ such
that $\fra^r$ is not contained in the $e^{\rm th}$ Frobenius power
of $J$
$$J^{[p^e]}:=(u^{p^e}\mid u\in J).$$
It is easy to see that $\nu^J(p^e)/p^e\leq\nu^J(p^{e+1})/p^{e+1}$,
and the $F$-\emph{threshold} of $\fra$ with respect to $J$ is
defined by
$$c^J(\fra):=\sup_{e}\frac{\nu^J(p^e)}{p^e}.$$
It is shown in \cite{MTW} that the set of $F$-thresholds of $\fra$
(with respect to various $J$) is precisely the set of jumping
numbers for the test ideals of $(X,Y)$. Note that the smallest
$F$-threshold is obtained for $J=\frakm$: this is an analogue of the
log canonical threshold that was introduced and studied by Takagi
and Watanabe in \cite{TW}.

\bigskip

There are several interesting results and questions relating the
invariants in characteristic zero and those obtained via reduction
mod $p$. To keep the notation simple we will work in the following
setup. Suppose that $\fra$ is an ideal in $A[x_1,\ldots,x_n]$, where
$A$ is the localization of $\ZZ$ at some integer. Let $Y$ be the
subscheme of $X=\AAA_A^n$ defined by $\fra$. If $p$ is a prime that
is large enough, then by reducing mod $p$ and localizing at
$(x_1,\ldots,x_n)$ we get a closed subscheme $Y_p$ in $X_p={\rm
Spec}\,{\mathbb F}_p[x_1,\ldots,x_n]_{(x_1,\ldots,x_n)}$ defined by
the ideal $\fra_p$.

The multiplier ideals of the pair $(X,Y)$ (more precisely, of its
extension to $\CC$) can be computed by a log resolution defined over
$\QQ$. After suitably localizing $A$ we may assume that the
multiplier ideals are defined over $A$, too. The following results
relate the reduction mod $p$ of the multiplier ideals with the test
ideals. They are due to Hara and Yoshida \cite{HY}, based on
previous work of Hara, Smith, Takagi and Watanabe.

\vfill\eject

\begin{theorem}\label{charp1}
With the above notation, if $p\gg 0$, then for every $\lambda$ we
have
$$\tau(X_p,\lambda\cdot Y_p)\subseteq \JJ(X,\lambda\cdot Y)_p.$$
\end{theorem}

Note that since our primes are large enough, the log resolution over
$\QQ$ induces by reduction mod $p$ log resolutions for $(X_p,Y_p)$.
The proof of Theorem~\ref{charp1} is based on the use of local
duality for the reduction mod $p$ of the log resolution. The proof
of the next result is more involved, using the approach of Deligne
and Illusie to the positive characteristic proof of the Kodaira
Vanishing Theorem.

\begin{theorem}\label{charp2}
With the above notation, for every $\lambda$ and for every $p\gg 0$
(depending on $\lambda$) we have
$$\tau(X_p,\lambda\cdot Y_p)=\JJ(X,\lambda\cdot Y)_p.$$
\end{theorem}

\smallskip

\begin{remark}
We reinterpret the above statements in terms of jumping numbers. For
simplicity, we restrict ourselves to the smallest such number: given
$\fra$ as above and $p\gg 0$, we want to compare the log canonical
threshold $c$ of the pair $(X,Y)$ in some small neighborhood of the
origin, with the $F$-pure threshold $c_p=c^{\frmm}(\fra_p)$.
Theorem~\ref{charp1} implies that for all $p\gg 0$ we have $c\geq
c_p$, while Theorem~\ref{charp2} implies that
$\lim_{p\to\infty}c_p=c$.
\end{remark}

\begin{example}\label{cusp}
Let $\fra$ be generated by $f=x^2+y^3$, whose log canonical
threshold is $\frac{5}{6}$. Let $p>3$ be a prime. One can show that
if
 $p\equiv 1$ (mod $3$), then the largest $r$ such that $f^r$
 does not lie in in $(x^{p^e}, y^{p^e})$ is given by $\nu(p^e)=\frac{5}{6}(p^e-1)$ for every
 $e\geq 1$, so that $c_p=\frac{5}{6}$. On the other hand, if
 $p\equiv 2$ (mod $3$), then
 $\nu(p)=\frac{5p-7}{6}$, while $\nu(p^e)=\frac{5p^e-p^{e-1}-6}{6}$
 for $e\geq 2$.
 Therefore in this case $c_p=\frac{5}{6}-\frac{1}{6p}$.
\end{example}

\begin{conjecture}
For every ideal $\fra$ in $A[x_1,\ldots,x_n]$ there are infinitely
many primes $p$ for which the $F$-pure threshold $c_p$ is equal to
the log canonical threshold $c$.
\end{conjecture}

For a discussion of this conjecture we refer to \cite{MTW}. We end
by mentioning a connection between the positive characteristic
invariants and the Bernstein-Sato polynomial. Suppose that $f\in
A[x_1,\ldots,x_n]$ is as above. We know that the Bernstein-Sato
polynomial $b_f(s)$ has rational roots, and in fact, one can show
that one can find an equation (\ref{bernstein1}) as in the
definition of $b_f(s)$ with $P$ having rational coefficients.
Therefore, after suitably localizing $A$, we may assume that both
$b_f$ and $P$ have coefficients in $A$ and that (\ref{bernstein1})
holds over $A$. It follows that if $p$ is a large enough prime, we
get a similar equation over ${\mathbb F}_p$.

Consider now an ideal $J$ in  the ring ${\mathbb
F}_p[x_1,\ldots,x_n]_{(x_1,\ldots,x_n)}$, such that $f_p$ lies in
the radical of $J$. Let us apply (\ref{bernstein1}) with
$s=\nu^J(p^e)$, the largest integer such that $f_p^r$ is not in
$J^{[p^e]}$. Since the ideal $J^{[p^e]}$ is a module over the ring
${\mathbb F}_p[x,\partial_x]$, we deduce that $b_f(\nu^J(p^e))\equiv
0$ (mod $p$). Therefore the functions $\nu^J$ give roots of $b_f$
mod $p$. Sometimes one can use this observation to find actual roots
of $b_f$.

\begin{example}
Let $f=x^2+y^3$. We have  described in Example~\ref{cusp} the
function $\nu=\nu^J$ when $J$ is the maximal ideal. If $p\equiv 1$
(mod 3), then $\nu(p^e)=\frac{5}{6}(p^e-1)$. The above discussion
implies that $p$ divides $b_f(-5/6)$. Since there are infinitely
many such primes, we deduce that $-\frac{5}{6}$ is a root of $b_f$.
Similarly, if $p\equiv 2$ (mod $3$), then it follows from the
formula for $\nu(p)$ that $-\frac{7}{6}$ is a root of $b_f$, and
from the formula for $\nu(p^e)$, with $e\geq 2$ that $-1$ is a root
of $b_f$. Therefore we have obtained all roots of the Bernstein-Sato
polynomial of $f$ by this method.
\end{example}

A similar picture can be seen in other examples, though at the
moment there is no general result in this direction. In \cite{BMS2}
this approach was used to describe all the roots of the
Bernstein-Sato polynomial of a monomial ideal. It would be very
interesting to find a more conceptual framework that would explain
the connection between the Bernstein polynomial and the invariants
in positive characteristic.

\smallskip

{\bf Acknowledgements}.
We would like to express a special thanks to Rob Lazarsfeld.
 Many results
in this paper were joint work with him.
Moreover, we have benefited from
many inspiring discussions.

\end{document}